\documentclass{article}
\usepackage{amsmath,amssymb,amsthm,dsfont,hyperref}
\begin{document}
	\title{On Fixed Points of a Map Defined by Alain Connes}
	\author{Hongbo Zhao}
	\date{}
	\maketitle
	\allowdisplaybreaks
	\newtheorem{theorem}{Theorem}
	\newtheorem{definition}{Definition}
	\newtheorem{example}{Example}
	\newtheorem{lemma}{Lemma}
	\newtheorem{proposition}{Proposition}
	\numberwithin{equation}{section}
	\section{Introduction}
	\par{
		In the theory of finite groups, the celebrated Feit-Thompson theorem states that any finite group of odd order is solvable \cite{FT63}. However, the proof is very long and technical. It is desirable to have a short proof.
	}
	\par{
		The Feit-Thompson theorem is equivalent to the fact that any non-trivial finite group $G$ of odd order has a non-trivial one-dimensional representation. Atiyah proposed the idea that one may simplify the proof of Feit and Thompson, using iteration processes to prove the existence of non-trivial one dimensional representations of $G$. In the paper \cite{C19}, Connes continued this idea and he defined a map $\Psi$ from the space of class functions $C(G)$ to $C(G)$ itself. The main observation is that one dimensional representations are fixed points of $\Psi$. One may expect to get a one-dimensional representation by iterating the map $\Psi$, and Connes made some analysis on iteration maps $\Psi^n$ when $n$ goes to infinity.
	}
	\par{
		The purpose of this paper, is to analyze the maps $\Lambda_t$ and $\Psi_t$ to be defined in Section 2 and Section 3, where $t$ is either a formal variable or a complex number. The map $\Psi$ is the same as $\Psi_t$ when $t$ equals $-1$, and we can show that one-dimensional representations of $G$ are also fixed points of $\Psi_t$ where $t$ is viewed as a formal variable. Our main consideration is that, instead of studying the fixed points of $\Psi$ which are not easy to handle, we consider the map $\Psi_t$  where $t$ is a formal variable, and we can characterize its fixed points explicitly. Moreover we can compare the fixed points with one dimensional representations for some concrete examples.  
	}
	\par{
		The main results of this paper is the following. In Section 1, we introduce the maps $\Lambda_t$ and $\Psi_t$ by analogy. Our first result, Theorem 1, shows that when $t$ is viewed as a formal variable, the fixed points of $\Psi_t$, except the zero solution, is a finite abelian group denoted as $A(G)$. In Section 2 we view $t$ as complex number at the moment and study the analytic properties of $\Lambda_t$ and $\Psi_t$. This is related to Section 3 of \cite{C19}, and we prove Theorem 2 which is about analytic continuations of $\Lambda_t$ and $\Psi_t$. In Section 3 we calculate $A(G)$ for $G=S_n$, $A_n$ and the results are given by Theorem 3 and Theorem 4. We also calculate $A(G)$ for $G$ being some types of finite abelian groups. From the illustrating examples, the calculation shows that there are many fixed points of $\Psi_t$ which are not one dimensional representations.
	}
	\section{The Maps $\Lambda_t$ and $\Psi_t$}
	\par{
		In this section we will prove the first main result, Theorem 1. We introduce and study the maps $\Lambda_t$ and $\Psi_t$ where $t$ is viewed as a formal variable in this section. To define the maps we give a brief overview of the main idea of the paper \cite{C19} by Connes, and recall some `natural operations' such as $\lambda$-operations and Adams operations. Finally we prove Theorem 1 which gives the algebraic characterization of the fixed points of the map $\Psi_t$. We remark that the $\lambda$-operations and Adams operations used in this paper are just a particular case of general $\lambda$-operations and Adams operations in $K$-theory.
	}
	\par{
		 Let $G$ be a finite group and we assume that the ground field $k=\mathbb{C}$. We use $Rep(G)$ to denote the set of equivalence classes of finite dimensional $G$-modules, and $R(G)$ denotes the Grothendieck ring of $Rep(G)$. For a $G$-module $V$ we use $[V]$ to denote the coresponding equivalence class in $Rep(G)$. It is also well known that the complexified Grothendieck ring $R(G)_{\mathbb{C}}$ is isomorphic to $C(G)$, the ring of class functions on $G$.
	}
	\par{
		Recall that there is a $\lambda$-ring structure on $R(G)$ consists of a family of maps $\lambda^k:R(G)\rightarrow R(G),k\in\mathbb{Z}_{\geq 0}$ called $\lambda$-operations, which are given by the wedge powers:
		$$
		\lambda^k(v):=[\Lambda^k(V)]\text{ for all }v=[V]\in Rep(G).
		$$
		The $\lambda$-ring structure on $R(G)$ has the following properties: For any $v,w\in R(G)$,
		$$
		\lambda^0(v)=1,\,\lambda^1(v)=v,\,\lambda^n(1)=0\text{ for all }n>1,
		$$
		and we also check that
		\begin{align}
		\lambda^k(v+w)=\sum_{0\leq i\leq k}\lambda^i(v)\lambda^{k-i}(w).\label{lambda}
		\end{align}
		holds.
	}
	\par{
		Let $t$ be a formal variable and we define a map $\Lambda_t:R(G)\rightarrow R(G)[[t]]$ by
		$$
		\Lambda_t:=\sum_{k\geq 0}\lambda^kt^k,
		$$
		which uses a single formal power series representing the whole $\lambda$-ring structure. It follows from (\ref{lambda}) that $\Lambda_t(\cdot)$ satisfies:
		$$
		\Lambda_t(v+w)=\Lambda_t(v)\Lambda_t(w)\text{ for all }v,w\in R(G).
		$$
	}
	\par{
		In \cite{C19} Connes considered the case when $t$ is specialized to the number $-1$, and he uses the notation
		$$
		\Lambda:=\Lambda_{-1}.
		$$
		He also defined the map $\Psi$ by
		$$
		\Psi:=I-\Lambda
		$$
		where $I$ denotes the trivial representation. The main observation of Connes which motivates the definition of the map $\Psi$ is the following: Suppse $v=[V]\in R(G)$ is a one dimensional representation of $G$, then it is clear that 
		$$
		\Lambda^k(V)=0\text{ for all }k\geq 2,
		$$
		and therefore
		$$
		\Lambda(v)=I-v.
		$$
		That is, if $v$ is a one dimensional representation then $v$ is a fixed point of the map $\Psi$:
		$$
		\Psi(v)=v.
		$$
		
	}
	\par{
		The main idea of our consideration is the following: If $v$ is a one dimensional representation and we define
		\begin{align}
		\Psi_t:=\frac{1}{-t}(\Lambda_{-t}-I), \label{psi}
		\end{align}
		then
		$$
		\Lambda_{-t}(v)=I-vt
		$$
		and therefore $v$ is also a fixed point of $\Psi_t$:
		$$
		\Psi_t(v)=v .
		$$
		We note that here $-t$ can be divided because
		$$
		(\Lambda_{-t}(v)(x)-1)\in t\mathbb{C}[[t]]
		$$
		for any $v\in R(G)$ and $x\in G$. It is clear that being a fixed point of $\Psi_t$ is a much stronger condition than being a fixed point of $\Psi=\Psi_{-1}$, therefore It deserves to study the fixed points of $\Psi_t$ first. 
	}
	\par{
		We need to extend the definitions of $\Lambda_t$ and $\Psi_t$ from $R(G)$ to $C(G)$, i.e., the space of class functions . We note that $\Lambda_t$ is not an additive homomorphism from $R(G)$ to $R(G)[[t]]$. To solve this problem, first we need to use Adams operations, and second, we need to introduce the `formal logartithm wedge' $log(\Lambda_{-t})$, which is a linear map. 
	}
	\par{
		The Adams operations $\psi^k$ are defined by evaluating the traces of $x^k$: For any element $v=[V]\in R(G)$ and $k\in\mathbb{Z}_{\geq 0}$, the $k$-th Adams operation $\psi^k:R(G)\rightarrow R(G)$ is defined by the condition that
		\begin{align}
		\psi^k (v)(x)=Tr(x^k)_{|V} \label{trace}
		\end{align}
		holds for any $x\in G$. Moreover the Adams operations $\psi^k$ are extended to linear maps on $C(G)$ such that
		$$
		\psi^k(f)(x)=f(x^k)
		$$
		holds for any $k\in\mathbb{Z}_{\geq 0}$ and $x\in G$. We note that linearity requires that
		$$
		\psi(f+g)=\psi^k(f)+\psi^k(g),\,\psi^k(cf)=c\psi^k(f)
		$$
		hold for any $c\in \mathbb{C}$ and $f,g\in C(G)$.
		It is clear that the Adams operations satisfy
		$$
		\psi^1(v)=v,\,\psi^k\psi^l=\psi^{kl},\,\psi^0(v)=dim(v).
		$$
		for any $v=[V]$ with $V$ a $G$-module.
	}
	\par{
		The Adams operations are related to the $\lambda$-operations as follows: For any $v\in R(G)$ we define
		\begin{align}
		\psi_t(v):=\sum_{k\geq 0}\psi^k(v)t^k=\psi^0(v)-t\frac{d}{dt}log(\Lambda_{-t})(v).\label{Adams}
		\end{align}
		Here (\ref{Adams}) is understood as an identity evaluated for any $x\in G$:
		$$
		\psi_t(v)(x)=\psi^0(v)-t\frac{d}{dt}log(\Lambda_{-t})(v)(x),
		$$
		and the logarithm $log(\Lambda_{-t})$ is defined in the formal sense: \begin{align}
		log(\Lambda_{-t})(v):=-\sum_{k\geq 1}\frac{t^k}{k}\psi^k(v).\label{log}
		\end{align}
		The identity (\ref{log}) is equivalent to the following where $\Lambda_{-t}(x)$ is expressed in terms of Adams operations:
		$$
		\Lambda_{-t}(v)=exp(-\sum_{k\geq 1}\frac{t^k}{k}\psi^k(v)),
		$$
		which is also used as an equivalent definition of Adams operations by some people.
	}
	\par{ 
		Here we do not prove that the definitions using (\ref{Adams}) or (\ref{log}) are the same as the previous one using (\ref{trace}) and interested readers can consult \cite{A67} for general discussions. We just remark that they have the following elentary analogue: Suppose $B\in End(\mathbb{C}^n)$ is a linear map with eigenvalues $\lambda_1,\cdots,\lambda_n$, then the coefficients of the characteristic polynomial
		$$
		det(I_n-Bt)
		$$
		which are elementary symmetric functions of eigenvalues, can be expressed using 
		$$
		Tr(B^k)=\sum_i \lambda_i^k,\,k=1,\cdots,n,
		$$
		which are homogeneous symmetric functions of eigenvalues, as we note that
		\begin{align*}
		1-\lambda_i t=exp(log(1-\lambda t))=exp(-\sum_{k\geq 1}\frac{\lambda_i^k}{k}t^k)
		\end{align*}
		for a single eigenvalue. Therefore
		$$
		det(I_n-Bt)=\prod_i(1-\lambda_i t)=exp(-\sum_{k\geq 1,i=1,\cdots,n}\frac{\lambda_i^k}{k}t^k).
		$$
		The $\lambda$-operations are similar to elementary symmetric functions while Adams operations are similar to homogeneous symmetric functions.
	}
	\par{
		Using (\ref{log}) we can extend $log(\Lambda_{-t})(\cdot)$ to be a linear map from $C(G)$ to $C(G)[[t]]$ because the right hand side is a sum of Adams operations $\psi^k$ which are linear, and therefore we can also extend $\Lambda_{-t}$ to be a map from $C(G)$ to $C(G)[[t]]$ through the exponential map:
		$$
		\Lambda_{-t}(v)(x):=exp(log(\Lambda_{-t})(v)(x))\,\text{for all }x\in G.
		$$
		It is also clear that by (\ref{psi}), the map $\Psi_t$ is also defined on $C(G)$.
	}
	
	\par{
		The fixed point of $\Psi_t$ has the following simple characterization.
		\begin{theorem}
			A class function $f\in C(G)$ is a fixed point of $\Psi_t$ if and only if $f$ satisfies:
			$$
			f(x^k)=f(x)^k
			$$
			for all $x\in G,\,k\in\mathbb{Z}_{\geq 0}$.
		\end{theorem}
		\textbf{Proof.} This can be shown by a direct calculation. For arbitrary $x\in G$, because
		$$
		-\sum_{k\geq 1}\frac{t^k}{k}\psi^k=log(\Lambda_{-t})=log(I-t\Psi_{t}),
		$$
		therefore the assumption that $f$ is a fixed point of $\Psi_t$ means that
		$$
		\Psi_{t} f(x)=f(x)
		$$
		holds for any $x\in G$, and it implies that
		\begin{align*}
		=&-\sum_{k\geq 1}\frac{t^k}{k}f(x^k)\\
		=&-\sum_{k\geq 1}\frac{t^k}{k}\psi^k (f)(x)= log(I-t\Psi_{t} (f))(x)\\
		=&log(1-tf(x))=-\sum_{k\geq 1}\frac{t^k}{k}f(x)^k,
		\end{align*}
		in which we note that we have the following formal identity:
		$$
			log(1-u)=-\sum_{k\geq 1}\frac{u^k}{k},
		$$
		and for the trivial representation $I$,
		$$
			I(x)=1
		$$
		holds for any $x\in G$. The proof is concluded by comparing the coefficients. 
	}
	\section{Analytic Properties of $\Psi_t$}
	\par{
		In this section we discuss analytical aspects of the maps $\Lambda_{-t}$ and $\Psi_t$, that is,  $t$ is viewed as a complex variable and therefore $\Lambda_{-t}$ and $\Psi_t$ will be analytic maps from $C(G)$ to $C(G)$ parametrized by $t$.}
	\par{
		 The discussion of this part follows Section 3 of Connes's paper \cite{C19},  which reviews some results of Euler on divergent series. Our main idea is to do analytic continuations. More precisely for fixed $f\in C(G)$ and $x\in G$ we view $\Lambda_{-t}(f)(x)$ and $\Psi_t(f)(x)$ as complex analytic functions of $t$ and we do analytic continuations starting from the obvious convergence domain. Because $\Psi_t$ is expressed in terms of $\Lambda_{-t}$ and therefore we only consider $\Lambda_{-t}$. The main result of this section is Theorem 2 and we will give the proof.
	}
	\par{
		Recall that that for $f\in C(G)$ and $x\in G$,
		$$
		log(\Lambda_{-t}(f))(x)=-\sum_{k\geq 1}\frac{t^k}{k}\psi^k(f)(x)=-\sum_{k\geq 1}\frac{t^k}{k}f(x^k),
		$$
		and
		$$
		(\Lambda_{-t})(f)(x)=exp(-\sum_{k\geq 1}\frac{t^k}{k}f(x^k)).
		$$
		It is clear from the expression that $log(\Lambda_{-t}(f))(x)$ and $(\Lambda_{-t})(f)(x)$ are absolutely convergent on the open unit disk $|t|<1$, and we proceed by doing the analytic continuations.
	}
	\par{
		Let $n=|G|$ denote the order of $G$, then by taking the derivative with respect to $t$ we obtain
		\begin{align*}
		\frac{d}{dt}(\Lambda_{-t})(f)(x)=&(\Lambda_{-t})(f)(x)\cdot(-\frac{d}{dt}\sum_{k\geq 1}\frac{t^k}{k}f(x^k))\\
		=&(\Lambda_{-t})(f)(x)\cdot(-\sum_{r=0,\cdots,n-1,k\geq 0}t^{kn+r}f(x^{r+1}))\\
		=&(\Lambda_{-t})(f)(x)\cdot(-\sum_{r=0,\cdots,n-1}f(x^{r+1})\frac{t^r}{1-t^n}).
		\end{align*}
		on $|t|<1$. This means that $(\Lambda_{-t})(f)(x)$ satisfies the linear ordinary differential equation
		$$
		\frac{d}{dt}(\Lambda_{-t})(f)(x)=(-\sum_{r=0,\cdots,n-1}f(x^{r+1})\frac{t^r}{1-t^n})\cdot (\Lambda_{-t})(f)(x),
		$$
		or equivalently
		\begin{align}
		\frac{d}{dt}(log(\Lambda_{-t}))(f)(x)=(-\sum_{r=0,\cdots,n-1}f(x^{r+1})\frac{t^r}{1-t^n}).\label{residue}
		\end{align}
		Let $\omega_n$ be a $n$-th primitive root of unity, then we note that the function $\sum_{r=0,\cdots,n-1}f(x^{r+1})\frac{t^r}{1-t^n}$ 
		has only simple poles at roots of unity $t=\omega^p_n$, $p=1,\cdots,n$ on the complex plane. Therefore by the analytic theory of ordinary differential equations, $(\Lambda_{-t})(f)(x)$ can be analytically continued to a multi-valued analytic function of $t$ on the complex plane with the $n$-th roots of unity being regular singularities. 
	}
	\par{
		To be more precise we need to choose a basis of $C(G)$ to give the explicit form of the linear map $log(\Lambda_{-t})$. We note that elements in the same conjugacy class are sent to the same conjugacy class under the power map $x\mapsto x^k$, therefore it is natural to choose the characteristic functions on the conjugacy classes as a basis. Suppose $C_1,\cdots,C_l$ are conjugacy classes of $G$, $c_i\in C_i,\,i=1,\cdots,l$ are representatives of $C_i$ and $[c_i]:=C_i$. For $r=0,\cdots n-1$, $j=1,\cdots l$, let $d_j^r$ be the number labelling the $r$-th power of the conjugacy class $C_j$, that is,
		$$
		[c_j^r]=[c_{d_j^r}].
		$$
		Let $\chi_i,i=1,\cdots,l$ denote the characteristic function of the class $C_i$:
		$$
		\chi_i(C_j)=\delta_{i,j}.
		$$
		We also introduce the following `two-variable Hurwitz like' function $H(\alpha,u)$ which is a variant of the Hurwitz function mentioned in \cite{C19}:
		$$
		H(\alpha,u):=\sum_{k\geq 0}\frac{u^{k}}{k+\alpha}.
		$$
		Then
		\begin{align*}
		&log(\Lambda_{-t}(\chi_i))(c_j)\\
		=&-\sum_{k\geq 1}\frac{t^k}{k}\chi_i(c_j^k)\\
		=&-\sum_{r=0\cdots,n-1}\frac{\chi_i(c_j^{r+1})t^{r+1}}{n}H(\frac{r+1}{n},t^n)\\
		=&-\sum_{r=0\cdots,n-1}\frac{\delta_{i,d_j^{r+1}}t^{r+1}}{n}H(\frac{r+1}{n},t^n),
		\end{align*}
		which implies that
		$$
		log(\Lambda_{-t}(\chi_i))=-\sum_{r=0\cdots,n-1;j=1,\cdots,l}\frac{\delta_{i,d_j^{r+1}}}{n}H(\frac{r+1}{n},t^n)t^{r+1}\chi_{j}.
		$$
		In general if $f$ is expressed by
		$$
		f=f_1\chi_1+\cdots+f_l\chi_l,
		$$
		then 
		$$
		log(\Lambda_{-t}(f))=-\sum_{r=0\cdots,n-1;j=1,\cdots,l}\frac{f_{d_j^{r+1}}}{n}H(\frac{r+1}{n},t^n)t^{r+1}\chi_{j}.
		$$
		It is calculated that
		$$
		Res_{t=\omega_n^p}(\frac{d}{dt}log(\Lambda_{-t}(f))(c_j))=-\sum_{r=0\cdots,n-1}\frac{\omega_n^{pr}f_{d_j^{r+1}}}{n}.
		$$
		by (\ref{residue}). The function $\Lambda_{-t}(f))(c_j)$ will be a single-valued function around $\omega_n^p$ exactly when the above residue is an integer. In summary we have
		\begin{theorem}
			For any fixed $f\in C(G)$ and a conjugacy class $c$ of $G$ the function $\Lambda_{-t}(f)(c)$ can be analytically continued to a multi-valued analytic function with $n$-th roots of unity $1,\omega_n\cdots,\omega_n^{n-1}$ being regular singularities. For $f=f_1\chi_1+\cdots+f_j\chi_l$ $j=1,\cdots,l$ and $p=0,\cdots,n-1$, if the condition
			$$	
			k_p:=-\sum_{r=0\cdots,n-1}\frac{\omega_n^{pr}f_{d_j^{r+1}}}{n}\in\mathbb{Z}
			$$
			is satisfied, then $\Lambda_t(f)(c_j)$ can be analytically continued to a single valued function around $\omega_n^p$. Moreover the isolated singular point $t=\omega_n^p$ is removable if $k_p\in\mathbb{Z}_{\geq 0}$ and in this case it is a $k_p$-th order zero, and otherwise $t=\omega_n^p$ is a $-k_p$-th order pole.
		\end{theorem}
	}
	\par{
		In particular for $t=-1$, $\Lambda_{1}(f)(x)$ is well defined for arbitrary $f$ and $x$ if and only if $n$ is odd. This is exactly Proposition 2.1 in \cite{C19} which was proved in a different way. 
	}
	\par{
		We add a final remark to this section that, however, we are unable to characterize the fixed points of $\Psi_t$ where $t$ is a given number. It is clear that if $f\in C(G)$ is a fixed point of $\Psi_t$ where $t$ is a formal variable, then $f$ is automatically a fixed point of $\Psi_t$ for a specific $t$ if $\Psi_t$ is defined, but the converse is not true. Therefore we may expect that there are many more fixed points of $\Psi_t$ where $t$ is a fixed complex number. Based on some analysis in Section 7 of \cite{C19} on the Lambert functions, we conjecture that for `generic' $t$ the map $\Psi_t$ has infinitely many fixed points.
	}
	\section{The Fucntor $A(\cdot)$ and Calculations of Some Examples} 
	\par{
		The main content of this section is to define the abelian group $A(G)$ for an arbitrary finite group $G$, and we compute some examples of $A(G)$. 
	}
	\par{
		In the Section 3 we showed that a class function $f\in C(G)$ is a fixed point of $\Psi_t$ where $t$ is a formal variable, if and only if it is a homomorphism when restricted to any cyclic subgroups:
		$$
		f(x^n)=f(x)^n\text{ for all }x\in G,\,n\in\mathbb{Z}_{\geq 0}.
		$$
		Because $f(1)=f(1\cdot 1)=f(1)^2$, therefore $f(1)=1$ or $f(1)=0$. If $f(1)=0$ then $f(x)=0$ for all $x\in G$. In later discussions we only consider the case $f(1)=1$, and define
		$$
		A(G):=\{f|\,f\in C(G),f\neq 0,f(x^n)=f(x)^n\text{ for all }x\in G,\,n\in\mathbb{Z}\}.
		$$
		It is clear that $A(G)$ is an abelian group with the identity being the constant function $1_G$. The obvious (abelian) multiplication and invserse are
		$$
		(f\cdot g)(x)=f(x)g(x),\,f^{-1}(x)=\frac{1}{f(x)}\text{ for all }x\in G,\,f,g\in A(G).
		$$
		Moreover $A(G)$ is a fintie abelian group because the value of $f$ must be a $n$-th root of unity.
	}
	
	\par{
		We can show that $A(\cdot)$ satisfies some functorial properties. Suppose 
		$$
		F:G\rightarrow G'
		$$
		is a group homomorphism, then it induces a homomorphism of abelian group 
		$$
		A(F):A(G')\rightarrow A(G)
		$$
		by the pull back
		$$
		(A(F)(g))(x):=g(F(x))\text{ for any }g\in A(G'),x\in G,
		$$
		and it is checked that 
		$$
		g(F(\cdot))\in A(G)
		$$
		holds. Therefore from the categorical point of view we can view $A(\cdot)$ abstractly as a contravariant functor from the category of fintie groups to the category of finite abelian group.
	}
	\par{
		The group $A(G)$ has the following property which will be useful in some computations.
		\begin{proposition}
			if $G,H$ are two fintie groups with
			$$
			(|G|,|H|)=1,
			$$
			then
			$$
			A(G\times H)\simeq A(G)\times A(H).
			$$
		\end{proposition}
		\textbf{Proof.} For any $a\in A(G\times H)$, we define
		$$
		a_1(g):=a((g,1)),\,a_2(h):=a((1,h)).
		$$
		Now we show that
		$$
		a((g,h))=a_1(g)a_2(h)
		$$
		which implies Proposition 1. Because $ord(g)||G|$, $ord(h)||H|$ and $(|G|,|H|)=1$, therefore there are integers $p,q$ such that
		$$
		p|G|+q|H|=1.
		$$
		Hence
		$$
		a((g,h))^{q|H|}=a(g^{q|H|},h^{q|H|})=a((g,1))=a_1(g),
		$$
		and similarly
		$$
		a((g,h)^{p|G|})=a_2(h).
		$$
		It is obvious that
		$$
		a((g,h))=a((g,h))^{p|G|}a((g,h))^{q|H|}=a_1(g)a_2(h)
		$$
		holds.
	}
	\par{
		Now we calculate $A(G)$ for some examples of finite groups $G$. First we consider the case $G=S_n$.
	}
	\par{
		\begin{theorem}
			The abelian group $A(S_n)$ is isomorphic to $(\mathbb{Z}/2\mathbb{Z})^{\#\mathcal{P}^{odd}_{2^*}(n)}$.
		\end{theorem}
		Here $\mathcal{P}^{odd}_{2^*}(n)$ is the set of certain restricted integer partitions of $n$ which will be introduced during the proof. As a preparation we analyze the power map on the set of conjugacy classes of $S_n$ first.
	}
	\par{
		It is well known that conjugacy classes of $S_n$ are classified by the corresponding cycle shapes of elements, which are labelled by (unodered) integer partitions of $n$. In later discussions we will identify a conjugacy class in $S_n$ with the corresponding integer partition, and for an element $x\in G$ we use $[x]$ to denote the corresponding conjugacy class or integer partition according to the context. We also use $ord(x)$ to denote the order of an element $x\in G$, and because elements in the same conjugacy class $C$ have the same order, we may say a conjugacy class or the corresponding partition $[x]$ has order $ord(x)$.
	}
	\par{
		It is obvious that disjoint cycles are still disjoint under the power map, so we only need to consider the power map of a single cycle. The following lemma is done by an elementary calculation.
		\begin{lemma}
			The $r$-th power of a cycle of length $l$ equals $(r,l)$-number of disjoint cycles of length all equal to $l/(r,l)$.  
		\end{lemma}
		As a consequence we have
		\begin{lemma}
			If $(ord(x),r)=1$, then
			$$
			[x^r]=[x].
			$$
		\end{lemma}
	}
	\par{
		
		\begin{proposition}
			$A(S_n)$ ($n>1$) is isomorphic to $(\mathbb{Z}/2\mathbb{Z})^N$ for some positive integer $N$. 
		\end{proposition}
		\textbf{Proof.} For any $x\in S_n$ suppose $ord(x)=2^k(2m+1)$ where $m,k$ are non-negative integers. Then $ord(x^{2^k})=2m+1$ and $ord(x^{2m+1})=2^k$. By Lemma 2 we have $[x^{2^k}]=[(x^{2^k})^2]$, therefore for any $a\in A(S_n)$
		$$
		a(x^{2^k})^2=a((x^{2^k})^2)=a(x^{2^k})
		$$
		and 
		\begin{align}
		a(x)^{2^k}=a(x^{2^k})=1\label{id1}
		\end{align}
		hold. Similarly $[x^{2m+1}]=[(x^{2m+1})^3]$, hence
		$$
		a(x^{2m+1})^{3}=a((x^{2m+1}).
		$$
		and we have
		\begin{align}
		a(x)^{2(2m+1)}=1. \label{id2}
		\end{align}
		Combine (\ref{id1}) with (\ref{id2}) we have
		$$
		a(x)^2=1.
		$$
		This implies that any element in $A(S_n)$ is of order 2, so we conclude the proof by the classification of finite abelian groups.
	}
	\par{
		For convenience we identify $\mathbb{Z}/2\mathbb{Z}$ with the multiplicative group $\{1,-1\}$ and we proceed to determine the number $N$. We note that 
		$$
		a=a^2
		$$
		for any $a\in A(S_n)$ by Proposition 2, and for any $x\in S_n$ with $ord(x)=2^k(2m+1)$, $k,m$ non-negative integers we have
		$$
		a(x)=a(x)^{2m+1}=a(x^{2m+1}).
		$$
		Therefore any $a\in A(S_n)$ is determined by its value on the conjugacy classes whose orders are of the form $2^k$, $k\in\mathbb{Z}_{\geq 1}$ (the case $k=0$ is trivial). We note that the order of a partition is a positive integer power of $2$ if and only if the each number in the partition is a power of $2$, and at least one number is not equal to $1$. We also note that an odd power map induces the identity map on the set of conjugacy classes whose orders are powers of $2$, hence on this set we only need to consider the square map.
	}
	\par{
		We call $x\in S_n$ or the corresponding conjugacy class (or integer partition) $[x]$ a square if there exists $y$ with $x=y^2$. From the previous discussion we know that if $x=y^{2}$ for some $y$ then $a(x)=a(y)^2=1$ for all $a\in A(S_n)$. Therefore any $a\in A(S_n)$ is determined by the value on the set of non-square conjugacy classes with orders equal a positive power of $2$. 
	}
	\par{
		Let $\mathcal{P}_{2^*}(n)$ denote the set of integer partitions of $n$ whose orders are positive integer powers of $2$, or equivalently, any integer not equals $1$ and appears in the partition, is a positive integer power of $2$. We further use $\mathcal{P}^{odd}_{2^*}(n)$ to denote the subset of $\mathcal{P}_{2^*}(n)$ such that in the partition, there exists at least one integer $k\geq 1$ which occurs odd number of times. 
	}
	\par{
		We have the following lemma.
		\begin{lemma}
			For any $[x]\in \mathcal{P}_{2^*}(n)$, $x$ is not a square if and only if $[x]\in \mathcal{P}^{odd}_{2^*}(n)$.
		\end{lemma}
		This can be verified using Lemma 1.
	}
	\par{
		We give some examples explaining the above definitions and Lemma 3. Let $n=12$, the following partitions are some elements of $\mathcal{P}_{2^*}(12)$:
		\begin{align*}
		&P_1=1+1+1+1+2+2+4,\\
		&P_2=4+4+4,\\
		&P_3=4+8,\\
		&P_4=2+2+4+4,\\
		\end{align*}
		while
		\begin{align*}
		&P_4=1+2+2+3+4,\\
		&P_5=12.
		\end{align*}
		are not but we note that
		$$
		P_4^3=P_1,P_5^3=P_2.
		$$
		The partitions $P_1$, $P_2$ and $P_3$ are elements of $\mathcal{P}^{odd}_{2^*}(12)$ because for example, $4=2^2$ occurs one times in $P_1$ and $P_3$ and three times in $P_2$, while $P_4$ is not. Moreover $P_4$ is a square because
		$$
		P_4=P_3^2,
		$$
		while $P_1$, $P_2$ and $P_3$ are not squares. 
	}
	\par{
		\begin{proposition}
			Any map $a:\mathcal{P}^{odd}_{2^*}(n) \rightarrow \{1,-1\}$ can be extended to an element in $A(S_n)$.
		\end{proposition}
		From previous discussions we have shown that any $a\in A(S_n)$ is determined by its values on $\mathcal{P}^{odd}_{2^*}(n)$. Therefore Theorem 2 is a corollary of Proposition 3 as $N=\#\mathcal{P}^{odd}_{2^*}(n)$.
	}
	\par{
		\textbf{Proof of Proposition 3.} It is obvious that if $a\in A(S_n)$ then we should have
		$$
		a(x)=1
		$$
		for all $[x]\in \mathcal{P}_{2^*}(n)/\mathcal{P}_{2^*}^{odd}(n)$ because they are squares, and trivially we have
		$$
		a(1)=1.
		$$
		For a general element $x$ with $ord(x)=2^{k}(2m+1)$ we set
		$$
		a(x)=a(x^{2m+1})
		$$
		as either $x^{2m+1}=1$ or $[x^{2m+1}]\in \mathcal{P}_{2^*}(n)$. A direct case by case check shows that such $a$ is a class function which satisfies
		$$
		a(x)^n=a(x^n)
		$$
		for all $x\in G$ and $n\in\mathbb{Z}$ and therefore
		$$
		a\in A(S_n),
		$$
		which concludes the proof of Proposition 3 and Theorem 2.
	}
	\par{
		Next we do the calculation for the alternating group $G=A_n$. It is clear that for any conjugacy class $C$ of $S_n$, either $C\subseteq A_n$ or $C\cap A_n$ is empty holds, and moreover $C\subseteq A_n$ if and only if the corresponding partition has even number of even integers. The following lemma characterizes the conjugacy classes of $A_n$, and a sketch of the proof can be found in some textbooks, for example, in Exercise 19-21, p.131 of \cite{DF04}.
		\begin{lemma}
			Any conjugacy class of $A_n$ is contained in a conjugacy class of $S_n$. More precisely, if $C$ is a conjugacy class of $S_n$ and $C\subseteq A_n$, then there are two possibilities
			\begin{enumerate}
				\item $C$ is a conjugacy class of $A_n$, and we call $C$ non-split.
				\item $C$ splits, that is, $C$ is a disjoint union of two conjugacy classes of $A_n$, which are denoted as $C_-$, $C_+$ by convention. For any $x\in C_-$ and $\tau\in S_n/A_n$, we have $\tau x\tau^{-1}\in C_+$ so that $C_-$ and $C_+$ have the same size.
			\end{enumerate} 
		\end{lemma}
	}
	\par{
		From now on we will fix the notation that $C$ denotes a conjugacy class of $S_n$ with $C\subseteq A_n$, and $C_-$, $C_+$ denotes the conjugacy classes of $A_n$ which $C$ splits into if $C$ splits. For any $x\in A_n\subseteq S_n$, we also use $[x]'$ to denote the corresponding conjugacy class in $A_n$.
	}
	\par{
		A more useful lemma using the integer partition type is the following 
		\begin{lemma}
			$C$ splits if and only if for an element $x\in C$ (and therefore for any $x\in C$ also), the centralizer $C_{S_n}(x)\subseteq A_n$. More precisely $C_{S_n}(x)\subseteq A_n$ if and only if any number in the partition $[x]$ is odd, and all numbers are distinct.
		\end{lemma}
		We remark that the last sentence `all numbers are distinct' includes the case that the partition $[x]$ has only one number.
	}
	\par{
		The additional work we need to do for $A_n$ is that we need to determine the power map by taking account into the split conjugacy classes. The following fact is useful in the calculations: For any cycle $(i_1,\cdots,i_l)\in S_n$ and $\tau\in S_n$ we have
		\begin{align}
		\tau\cdot(i_1,\cdots,i_l)\cdot \tau^{-1}=(\tau(i_1),\cdots,\tau(i_l)).\label{conj}
		\end{align}
	}
	\par{
		We note that if $C$ is non-split, then for any integer $p$ and $x\in C$, $[x^p]$ is also non-split by Lemma 1 and Lemma 5. In this situation the power map is determined by Lemma 1, and therefore we need to discuss the case when $C$ splits. With this assumption we note that $r:=ord(C)$ is odd.
		\par{
			We first consider the subcase that $p$ is an integer which is not coprime to $r$. We note that for any $x\in C$, $[x^p]$ is also non-split by Lemma 1 and Lemma 5, therefore in this subcase the power map can be also be determined using Lemma 1.
		}
		\par{
			We are left to consider the subcase that $(p,r)=1$. In this subcase $C=[x^p]$ for any $x\in C$ by Lemma 1 and Lemma 5, and we need to determine the induced $p$-th power map on the set $\{C_-,C_+\}$. We note that the map $x\mapsto x^p$ is a bijection on the set $C$ by our assumptions, therefore by Lemma 4 and Lemma 5, for any $x\in C_-$ the induced map on the two element set $\{C_-,C_+\}$ is either 
			$$
			[x^p]'=[x]_-,\,[x^p]'=[x]_+,
			$$
			or
			$$
			[x^p]'=[x]_+,\,[x^p]'=[x]_-.
			$$
			In the first case the induced map is the identity map, and in the second case we call the induced map an exchange map as it exchanges $C_-$ and $C_+$. Moreover, there is a homomorphism
			$$
			(\mathbb{Z}/r\mathbb{Z})^{\times}\rightarrow \mathbb{Z}/2\mathbb{Z},
			$$
			in which the left hand side is the group of $p$-th power map, $(p,r)=1$, and the right hand side is the order two group generated by the exchange map.
		}
		\par{
			In principle we can determine the type of the induced map as follows. For the particular case when $C$ is represented by a single cycle $x$ of length $2m+1$, $m\in\mathbb{Z}_{\geq 1}$ and therefore $C$ splits:
			$$
			x=(1,\cdots,2m+1)\in C.
			$$ 
			Then the $p$-th power $x^p$ is also a single cycle of length $2m+1$, and we can find $\tau$ such that
			$$
			x^p=\tau x \tau^{-1}
			$$
			using (\ref{conj}). Therefore the map induced by $x\mapsto x^p$ is the identity map if $sgn(\tau)=1$, and otherwise the exchange map if $sgn(\tau)=-1$, by Lemma 4 and Lemma 5.
		}
		\par{
			It is complicated to determine the induced map for general $p$ and $m$ using the above calculation directly, but we are able to determine the map for the particular case when $p=2$ for all $m$, and it is sufficent to use this determining $A(A_n)$.
			\begin{lemma}
				Suppose $x$ is a cycle of length $2m+1$ and therefore $[x]$ splits. Then the square map $x\mapsto x^2$ induces the identity map on $\{[x]_-,[x]_+\}$ if $m$ mod $4$ equals $0,3$, and otherwise induces the exchange map if $m$ mod $4$ equals $1,2$.
			\end{lemma}
			To show this we only need to note that
			$$
			(1,2,\cdots,m,m+1,m+2,\cdots,2m+1)^2=(1,3,\cdots,2m+1,2,4,\cdots,2m),
			$$
			and by (\ref{conj}) the two elements are conjugated by the permutation
			$$
			\tau=\begin{pmatrix}
			1&2&\cdots&m&m+1&m+2&\cdots&2m+1\\
			1&3&\cdots&2m+1&2&4&\cdots&2m
			\end{pmatrix}
			$$
			whose sign equals
			$$
			(-1)^{\frac{m(m+1)}{2}}.
			$$
		}
		\par{
			In general it is easy to establish the following
			\begin{proposition}
				Suppose $C$ splits and therefore by Lemma 5, $C$ has the integer partition type
				$$
				(2m_1+1)+(2m_2+1)+\cdots+(2m_s+1)
				$$
				where $m_i$ are all distinct non-negative integers. Let
				$$
				m=\sum\frac{m_i(m_i+1)}{2},
				$$
				then the square map $x\mapsto x^2$ induces the identity map if $m$ is even and otherwise it induces the exchange map if $m$ is odd. 
			\end{proposition}
		}
		\par{
			In particular we have
			\begin{lemma}
				If $ord(C)$ is odd, then the $4$-th power map induces the identity map on $C$, whenever $C$ splits or not.
			\end{lemma}
		}
		\par{
			As a consequence we have
			\begin{proposition}
				The group $A(A_n)$ is a product of $2$-cyclic and $3$-cyclic groups.
			\end{proposition}
			\textbf{Proof.} The proof is similar to the proof of Proposition 3 with a little bit change. Suppose $x\in A_n$, and $ord(x)=2^k(2m+1)$ where $m,k$ are non-negative integers. First we consider the case that $ord(x)$ is even, therefore $k\geq 1$, and as before we have $ord(x^{2m+1})=2^k$ and $ord(x^{2^k})=2m+1$. We note that in this case $[x^{2^k}]$ and $[(x^{2^k})^2]$ are both non-split, hence we have $[x^{2^k}]=[(x^{2^k})^2]$ as before, and therefore  
			$$
			a(x^{2^k})^2=a((x^{2^k})^2)=a(x^{2^k}),
			$$
			which implies that
			$$
			a(x)^{2^k}=1.
			$$
			Similarly $[x^{2m+1}]$ is also non-split and has even order, hence $[x^{2m+1}]=[(x^{2m+1})^3]$.
			Therefore
			$$
			a(x^{2m+1})^{3}=a((x^{2m+1}),
			$$
			which implies
			$$
			a(x)^{2(2m+1)}=1.
			$$
			In conclusion we have
			$$
			a(x)^2=1
			$$
			for $x\in A_n$ with $ord(x)$ even.
		}
		\par{
			For the case $k=0$, $ord(x)$ is odd, we first note that
			$$
			a(x)=1
			$$
			if $[x]$ is non-split. If $[x]$ splits we note that by Lemma 7, the $4$-th power map is the identity map, hence
			$$
			a(x)^4=	a(x^4)=a(x),
			$$
			which implies
			$$
			a(x)^{3}=1.
			$$
		}
		\par{
			We also have more restrictions when $ord(x)$ is odd.
			\begin{proposition}
				Suppose $ord(x)$ is odd and $p\neq 3$ is a prime with $p|ord(x)$. Then
				$$
				a(x)=1.
				$$ 
			\end{proposition}
			We only need to note that $[x^p]$ must be non-split and $ord(x^p)$ is odd, therefore 
			$$
			a(x)^p=1.
			$$
			By Proposition 3,
			$$
			a(x)^3=1,
			$$
			therefore we conclude the proof as we note that $(p,3)=1$.
		}
		\par{
			The previous discussions tell that a $3$-cyclic group may be a subgroup of $A(A_n)$. We note that a necessary condition for $\mathbb{Z}/3\mathbb{Z}$ being a subgroup of $A(A_n)$ is that there is an integer partition of $n$ whose order is of the form $3^k$, $k\in\mathbb{Z}_{\geq 0}$ by Proposition 6. We note that this is equivalent to that $n$ is a sum of distinct non-negative integer powers of $3$, and it is also clear that if such a partition exists then it is unique, using $3$-adic representations. Therefore $\mathbb{Z}/3\mathbb{Z}$ occurs as most once in $A(A_n)$. 
		}
		\par{
			Let $\bar{\mathcal{P}}_{2^*}^{odd}(n)$ denote the subset of $\mathcal{P}_{2^*}^{odd}(n)$ consists of (unordered) integer partitions of $n$ such that, there are even number of even integers, that is, we only consider the partitions corresponding to conjugacy classes in $A_n$. The result is given by the following theorem.
			\begin{theorem}
				$A(A_n)$ is isomorphic to $ (\mathbb{Z}/2\mathbb{Z})^{\#\bar{\mathcal{P}}_{2^*}^{odd}(n)}$, or $(\mathbb{Z}/2\mathbb{Z})^{\#\bar{P}_{2^*}^{odd}(n)}\times (\mathbb{Z}/3\mathbb{Z})$. Moreover, the later case happens if and only if the following two conditions for $n$ are satisfied:
				\begin{enumerate}
					\item $n$ is of the form
					\begin{align}
					n=3^{k_1}+\cdots+\cdots3^{k_s}\label{3adic}
					\end{align}
					where $k_i,i=1,\cdots,s$ are all distinct non-negative integers.
					\item 
					The summation of power indices
					$$
					\sum_i k_i
					$$	
					is an odd integer.
				\end{enumerate}
			\end{theorem}
			In particular, there are infinitely many $n$ such that $\mathbb{Z}/3\mathbb{Z}$ is a subgroup of $A(A_n)$.
		}
		\par{
			We first prove the case that condition 1 for $n$ is not satisfied. In this case by the proof of Proposition 5 and Proposition 6, we have
			$$
			a(x)^2=1
			$$
			for all $x\in A_n$ with $ord(x)$ even, and
			$$
			a(x)=1
			$$
			for all $x\in A_n$ with $ord(x)$ odd. Because for any $x\in A_n$ with $ord(x)$ even, the conjugacy class $[x]$ is non-split, therefore the remaining part of the proof for this case is completely the same as the proof for $S_n$ by restricting everything to $A_n$ and we will not repeat.
		}
		\par{
			Now we suppose that condition 1 for $n$ holds. Let $C^*$ denote the distinguished conjugacy class in $S_n$ corresponding to the integer partition (\ref{3adic}), then it is clear that
			$$
			ord(C^*)=3^{\max\{k_i\}}.
			$$
			We define
			\begin{align*}
			&A':=\{a|a\in A(A_n),a(x)=1\text{ for all }x\in C^*\},\\
			&A^*:=\{a|a\in A(A_n),a(x)=1\text{ for all }x\notin C^*\}.
			\end{align*}
			Then we have
			$$
			A(A_n)\simeq A'\times A^*
			$$
			because the set 
			$$
			\{x|\,x\notin C^*\}
			$$
			is closed under the power map, and $x^k\notin C^*$ for $x\in C^*$ if and only if $k$ is not coprime to $3$. Moreover
			$$
			A'\simeq (\mathbb{Z}/2\mathbb{Z})^{\#\bar{\mathcal{P}}_{2^*}^{odd}(n)}
			$$  
			by the analysis the same as before, and
			$$
			A^*\simeq \mathbb{Z}/3\mathbb{Z}\text{ or }\{1\},
			$$
			because for any $x\in C^*$, $[x^3]$ is non-split and $ord(x)$ is odd, hence
			$$
			a(C^*_+)^3=a(C^*_-)^3=a([x^3])=1.
			$$
			It remains to determine $A^*$ explicitly by studying the $p$-th power map on $C^*$ with $(p,3)=1$.
		}
		\par{
			For convenience we also identify the cyclic group $\mathbb{Z}/3\mathbb{Z}$ with the multiplicative group $\{1,\omega,\omega\}$, $\omega^3=1$. It is clear that $A^*\simeq \mathbb{Z}/3\mathbb{Z}$ if and only if the following two assignments
			$$
			a(C^*_-)=\omega,\, a(C^*_+)=\omega^2
			$$
			or
			$$
			a(C^*_-)=\omega^2,\, a(C^*_+)=\omega
			$$
			are compatible with power maps, that is, the condition 
			\begin{align}
			a([x]')^p=a([x^p]') \label{power}
			\end{align}
			holds for any positive integer $p$ with $(p,3)=1$.
		}
		\par{
			We first assume that condition 2 of Theorem 4 is not satisfied, that is,
			$$
			\sum k_i
			$$
			is even. In this situation, by Proposition 4 it follows that the map on $\{C^*_-,C^*_+\}$ induced by the square map $x\mapsto x^2$ is the identity map, as we note that
			\begin{align}
			\begin{cases}
			\frac{3^k-1}{2}\equiv 1(mod\,4),\text{ if }k\text{ is odd},\\
			\frac{3^k-1}{2}\equiv 0(mod\,4),\text{ if }k\text{ is even}.
			\end{cases}\label{kpower}
			\end{align}
			Hence the square map induces the identity map on $\{C^*_-,C^*_+\}$ and 
			$$
			a(x^2)=a(x)
			$$
			for any $x\in C^*$. Therefore
			$$
			A^*\simeq \{1\}
			$$
			in this case.
		}
		\par{
			It remains to show that if condition 2 is also satisfied, then (\ref{power}) holds. In this situation, we note that the map induced by the square map on $\{C^*_-,C^*_+\}$ is the exchange map, using Proposition 4 and (\ref{kpower}) again.
		}
		\par{
			We need the following lemma from elementary number theory which reduces the analyze of $p$-th power map to the square map. The proof can be found, for example in \cite{IR90}, Chapter 4, Theorem 2.
			\begin{lemma}
				For $k\geq 1$ the group $(\mathbb{Z}/3^k\mathbb{Z})^{\times}$ consists of invertible elements in the modular $3^k$ arithemtic, is a cyclic group with $2$ being a generator. 
			\end{lemma}
			Because $(p,3)=1$ is equivalent to that (the image of) $p$ is an element in $(\mathbb{Z}/3^k\mathbb{Z})^{\times}$, therefore by Lemma 8, for any $(p,3)=1$ and a fixed integer $k\geq 1$, there is a non-negative integer $l$ such that
			$$
			p\equiv 2^l(mod\,3^k).
			$$
			For our purpose it is necessary to set $k=\max\{k_i\}$, and it is sufficient to check that (\ref{power}) holds for any $p=2^l$, $l$ an arbitrary positive integer. This is verified directly by using Proposition 4 as we note that
			$$
			\begin{cases}
			2^l\equiv 2(mod\,3),\text{ if }l\text{ is odd},\\
			2^l\equiv 1(mod\,3),\text{ if }l\text{ is even}.
			\end{cases}
			$$
			Moreover we see that this induces a nontrivial homomorphism 
			$$
			(\mathbb{Z}/3^k\mathbb{Z})^{\times}\rightarrow \mathbb{Z}/2\mathbb{Z}.
			$$
			Hence we conclude the proof of Theorem 4.
		}
		\par{
			Finally we calculate $A(G)$ for some cases of finite abelian group. It is well known that any finite abelian group is isomorphic to a direct product of abelian $p$-groups, $p$ a prime. By Proposition 1, it is sufficient to calculate $A(G)$ where $G$ is an abelian $p$-group.
			\begin{proposition}
				For $G=(\mathbb{Z}/p\mathbb{Z})^m$, 
				$$
				A(G)\simeq (\mathbb{Z}/p\mathbb{Z})^{\frac{p^m-1}{p-1}}.
				$$
			\end{proposition}
			We only need to note that non-trivial cyclic subgroups of $G=(\mathbb{Z}/p\mathbb{Z})^m$ are all cyclic groups of order $p$, and $G$ is the union of all the cyclic groups of order $p$. We only need to count how many cyclic $p$-groups are there, and each cyclic groups of order $p$ will contributes a summand isomorphic to $\mathbb{Z}/p\mathbb{Z}$. Suppose there are $N$ distinct cyclic $p$ subgroups, then it is obvious that
			$$
			Np-(N-1)=p^m
			$$ 
			by counting the number and therefore
			$$
			N=\frac{p^m-1}{p-1}.
			$$
			It is also easy to see that this is the same as counting the number of $1$-dimensional subspace of the $m$-dimensional space over the fintie field $F_p$. It is unfortunate that we are not able to determine the result for general cases, as an explicit characterization of the cyclic subgroup lattice of a general abelian $p$-group is not available.
		}
		\par{
			We conclude this paper by comparing $A(G)$ with one dimensional irreducible representations of $G$ for previous examples. Recall that $A_n$ is a simple group if $n\geq 5$, and in this case $A_n$ has only trivial one-dimensional representations. It is also well known that $S_n$, $n\geq 2$ has only two one-dimensional representations, the trivial representation and the sign representation. For the case of finite abelian groups, the irreducible representations of are all one-dimensional and determined by the corresponding direct summand of cyclic groups. We see that for the examples we have calculated, even if we require that $t$ is a formal variable, there are many fixed points of $\Psi_t$, which are not one-dimensional representations, therefore the maps $\Lambda_t$, $\Psi_t$ are still not strong enough to find non-trivial one-dimensional representations. It is necessary to find some possible refinement of $\Psi_t$ which may give furthur progress on this approach.
		}
		
	\end{document}